# ON FALSE DISCOVERY CONTROL UNDER DEPENDENCE

By Wei Biao Wu

*University of Chicago*

A popular framework for false discovery control is the random effects model in which the null hypotheses are assumed to be independent. This paper generalizes the random effects model to a conditional dependence model which allows dependence between null hypotheses. The dependence can be useful to characterize the spatial structure of the null hypotheses. Asymptotic properties of false discovery proportions and numbers of rejected hypotheses are explored and a large-sample distributional theory is obtained.

**1. Introduction.** Since the seminal work of Benjamini and Hochberg (BH) [2], the paradigm of false discovery control has been widely used in multiple hypothesis testing problems and it is often more useful than the classical Bonferroni-type method. Suppose that we want to test $n$ hypotheses $H_i$, $1 \le i \le n$. Write $H_i = 0$ if the $i$th null hypothesis is true and $H_i = 1$ if otherwise. Let $V$ be the number of erroneously rejected null hypotheses which are actually true and let $R$ be the total number of rejected hypotheses. The false discovery proposition (FDP) is defined as

$$(1) \qquad FDP = \frac{V}{R \vee 1} \qquad \text{where } a \vee b = \max(a, b),$$

and the false discovery rate (FDR) is defined as the expected value $\mathbb{E}(FDP)$.

We now briefly describe the BH procedure. Let $X_i$ be the marginal $p$-value of the $i$th test, $1 \le i \le n$, and let $X_{(1)} \le \cdots \le X_{(n)}$ be the order statistics of $X_1, \ldots, X_n$. Given a control level $\alpha \in (0,1)$, let

$$(2) \qquad R = \max\{i \in \{0, 1, \ldots, n+1\} : X_{(i)} \le \alpha i/n\},$$

where $X_{(0)} = 0$ and $X_{(n+1)} = 1$. The BH procedure rejects all hypotheses for which $X_{(i)} \le X_{(R)}$. If $R = 0$, then all hypotheses are accepted. Assume that $X_i$, $1 \leq i \leq n$, are independent and the $p$-value distribution is continuous; BH [2] proved that, if there are $N_0$ true null hypotheses, then $\mathbb{E}[V/(R \vee 1)] = \alpha N_0/n$. A popular framework for the false discovery control is the random effects model or the two-component mixture model (McLachlan and Peel [12]) in which the null hypotheses $H_i, 1 \leq i \leq n$, are assumed to be independent Bernoulli random variables. In particular, one assumes that $(X_i, H_i)$ are independent and identically distributed (i.i.d.) with

(3) $\quad \mathbb{P}(X_i \leq x|H_i = 0) = x, \qquad \mathbb{P}(X_i \leq x|H_i = 1) = G(x), \qquad 0 \leq x \leq 1,$

and $H_i \sim Bernoulli(\pi_1)$ [viz., $\mathbb{P}(H_i = 1) = \pi_1$ and $\pi_0 = 1 - \pi_1 = \mathbb{P}(H_i = 0)$]. Here $G$ is the distribution function of the $p$-value $X_i$ under alternative hypotheses. It is commonly assumed that $X_i \sim \text{uniform}(0, 1)$ if $H_i = 0$.

Due to the independence assumption, the classical random effects model or the two-component mixture model does not allow one to model spatial or location structures of the null hypotheses. In certain applications one expects that false null hypotheses occur in clumps, which are spatially clustered. In this case it is reasonable to expect that, if $H_i = 1$, then the nearby hypotheses $H_j$, where $j$ is close to $i$, are more likely to be false. In the negative dependence case the occurrence of $H_i = 1$ prevents nearby hypotheses from being false. Recently the multiple testing problem under spatial dependence has been considered by Qiu et al. [16] for microarray data and by de Castro and Singer [6] for geographical data.

In this paper we shall consider the problem of false discovery control without the independence assumption. In particular, we propose the conditional independence model: Let $(H_i)$ be a 0/1-valued stationary process, and, given $(H_i)_{i=1}^n$, $X_i$ are independent. The dependence is imposed on the hypotheses $(H_i)$. A simple relaxation of the independence assumption on $(H_i)$ is to impose a Markovian structure. In this case it is interestingly related to hidden Markov models (see Section 3).

As demonstrated in Storey, Taylor and Siegmund [19], Genovese and Wasserman [9], Chi [5] and Meinshausen and Rice [13] among others, the theory of empirical processes plays a useful role in the study of false discovery control. Recently Wu [24] considered empirical distribution functions for a wide class of stationary processes. In this paper we shall deal with the $p$-values arising from the aforementioned conditional independence model. In particular, we shall prove the validity of the BH procedure and present a distributional theory for $R$, the number of rejected hypotheses. We shall also establish a Bahadur-type asymptotic expansion for the false discovery proportion $V/(R \vee 1)$ and the weak convergence of false discovery processes to Gaussian processes.

The rest of the paper is structured as follows. Our dependence structure and main results are presented in Section 2 and proved in Section 4. Applications to Markov models and linear processes are given in Section 3.

**2. Main results.** We assume that $(H_s)_{s \in \mathbb{Z}^d}$ is a stationary random field and, for presentational simplicity, we shall consider testing hypotheses $H_s$ over $d$-dimensional cubes (cf. Condition 1). Results obtained in the paper can be generalized without essential difficulties to other types of regions. For a random variable $\xi$ write $\|\xi\| = \{\mathbb{E}(|\xi|^2)\}^{1/2}$. Denote by $\Rightarrow$ the weak convergence and by $N(\mu, \sigma^2)$ a normal distribution with mean $\mu$ and variance $\sigma^2$. Let $\mathcal{N}$ denote a standard normal random variable.

CONDITION 1. Let $(H_s)_{s \in \mathbb{Z}^d}$ be a stationary, $0/1$-valued random field. For $n_1, \ldots, n_d \in \mathbb{N}$ let the $d$-dimensional cube $C = \{1, 2, \ldots, n_1\} \times \cdots \times \{1, 2, \ldots, n_d\}$ and $n = n_1 n_2 \cdots n_d$. Write the sum $N_C = \sum_{s \in C} H_s$. Let $\pi_0 = \mathbb{P}(H_s = 0)$ and $\pi_1 = 1 - \pi_0$. Assume that, as $\min_{k \leq d} n_k \to \infty$, $\|N_C - n\pi_1\| = O(\sqrt{n})$, and the central limit theorem (CLT) $n^{-1/2}(N_C - n\pi_1) \Rightarrow N(0, \sigma^2)$ holds for some $\sigma^2 < \infty$.

In Section 3 we will present examples that Condition 1 is satisfied. With a slight abuse of notation, we write $(H_s)_{s \in C}$ as $(H_i)_{i=1}^n$, where $i = 1, \ldots, n$ corresponds to the lexicographic ordering of $s \in C$.

Under the conditional independence model we can have the representation

$$X_i = (1 - H_i)U_i + H_i G^{-1}(U_i), \tag{4}$$

where $U_i$ are independent and identically distributed (i.i.d.) uniform$(0, 1)$ random variables which are also independent of $(H_i)_{i=1}^n$, and $G^{-1}(u) = \inf\{x \in [0, 1] : G(x) \geq u\}$ is the inverse of $G$. Clearly (4) implies that the conditional distribution $[X_i | H_i = 0]$ is uniform$(0, 1)$ and $[X_i | H_i = 1]$ is $G$. If $(H_i)$ are independent, then (4) reduces to the random effects model. Our dependence paradigm is different from earlier ones adopted in Farcomeni [7] and Benjamini and Yekutieli [3].

Following Genovese and Wasserman [9], we consider the false discovery process

$$\Gamma_n(t) = \frac{n\Lambda_n(t)}{nF_n(t) + \prod_{i=1}^n \mathbf{1}_{X_i > t}}, \qquad 0 \leq t \leq 1, \tag{5}$$

where

$$\Lambda_n(t) = \frac{1}{n}\sum_{i=1}^n (1 - H_i)\mathbf{1}_{X_i \leq t} \quad \text{and} \quad F_n(t) = \frac{1}{n}\sum_{i=1}^n \mathbf{1}_{X_i \leq t}. \tag{6}$$

Then $F_n$ is the empirical process of $X_1, \ldots, X_n$ and $\Lambda_n$ can be interpreted as a marked empirical process. Let

$$\Delta_n(t) = \frac{1}{n}\sum_{i=1}^n H_i \mathbf{1}_{X_i \leq t} = F_n(t) - \Lambda_n(t). \tag{7}$$

Note that $H_i \mathbf{1}_{X_i \leq t} = H_i \mathbf{1}_{G^{-1}(U_i) \leq t}$ and $(1 - H_i) \mathbf{1}_{X_i \leq t} = (1 - H_i) \mathbf{1}_{U_i \leq t}$. Under the conditional independence model (4), we have for $0 \leq t \leq 1$ that

$$\Lambda(t) := \mathbb{E}\Lambda_n(t) = t\pi_0 \quad \text{and} \quad \Delta(t) := \mathbb{E}\Delta_n(t) = G(t)\pi_1.$$

To obtain large-sample properties of the false discovery process $\Gamma_n$, we need to establish an asymptotic theory for $\Lambda_n(t) - \Lambda(t)$ and $\Delta_n(t) - \Delta(t)$. Theorem 1 below concerns the weak convergence of $\sqrt{n}[\Lambda_n(t) - \Lambda(t)]$ and $\sqrt{n}[\Delta_n(t) - \Delta(t)]$ in a functional space. Let $\mathcal{D}[0,1]$ be the collection of functions which are right continuous and have left limits; let $\mathcal{D}^2[0,1] = \{(f_1, f_2) : f_1, f_2 \in \mathcal{D}[0,1]\}$. Assume throughout the paper that $G$ has a bounded density $g = G'$, namely, $\sup_{x \in [0,1]} g(x) < \infty$. Asymptotic results in Theorems 1–4 below are meant as $\min_{k \leq d} n_k \to \infty$.

THEOREM 1. *Assume Condition 1. Then there exist tight centered Gaussian processes $W_\Lambda(t)$ and $W_\Delta(t)$, $0 \leq t \leq 1$, such that the weak convergence*

$$(8) \quad (\sqrt{n}\{\Lambda_n(t) - \Lambda(t)\}, \sqrt{n}\{\Delta_n(t) - \Delta(t)\}) \Rightarrow (W_\Lambda(t), W_\Delta(t))$$

*holds in the space $\mathcal{D}^2[0,1]$.*

Since $F_n$ is a nondecreasing function and $X_{(j)}$ is the $j$th quantile of $F_n$, the value $R$ defined in (2) satisfies $R = \max\{0 \leq j \leq n : j/n \leq F_n(\alpha j/n)\}$. Let

$$(9) \quad \begin{aligned} \nu_{\text{BH}} &= \sup\{t \in [0,1] : t/\alpha \leq F_n(t)\} \quad \text{and} \\ \nu_0 &= \sup\{t \in [0,1] : t/\alpha \leq F(t)\}. \end{aligned}$$

It is easily seen that $R \leq n\nu_{\text{BH}}/\alpha < R + 1$. Let $f(x) = F'(x)$ and

$$(10) \quad \alpha_* = \frac{1}{f(0)} = \frac{1}{F'(0)} = \frac{1}{\pi_0 + \pi_1 g(0)}.$$

If $\pi_1$ and $g(0)$ are large, then $\alpha_*$ is small. Theorem 2 below describes asymptotic behavior of $\nu_{\text{BH}}$ and suggests a dichotomous phenomenon. It gives a Bahadur representation of $\nu_{\text{BH}}$ when $\alpha > \alpha_*$ and $R = O_\mathbb{P}(1)$ when $\alpha < \alpha_*$. At the boundary case $\alpha = \alpha_*$ we have an interesting nonstandard limiting distribution with a cubic root normalizing constant. In the case of random effects model in which $H_i$ are i.i.d., Chi [5] obtained interesting results on strong convergence properties of $R$ for the two cases $\alpha > \alpha_*$ and $\alpha = \alpha_*$. Chi also obtained a distributional result for $R$ when $\alpha < \alpha_*$ and argued that the number of rejected hypotheses is bounded even if there is a positive proportion of untrue null hypotheses. Chi's work shows the criticality phenomenon of false discovery rate controlling procedures.

THEOREM 2. *Assume Condition 1.*

(i) If $\alpha_*^{-1} > \alpha^{-1} > f(\nu_0)$, then

$$\nu_{\rm BH} - \nu_0 = \frac{F_n(\nu_0) - \nu_0/\alpha}{\alpha^{-1} - f(\nu_0)} + O_{\mathbb{P}}(n^{-3/4}). \tag{11}$$

Consequently $\sqrt{n}(\nu_{\rm BH} - \nu_0) \Rightarrow N(0, \sigma^2)$ for some $\sigma^2 < \infty$.

(ii) If $\alpha < \alpha_*$, then $R = O_{\mathbb{P}}(1)$.

(iii) If $\alpha = \alpha_*$ and $c_0 = -f'(0)/[2\sqrt{f(0)}] > 0$, then

$$n^{1/3}\nu_{\rm BH} \Rightarrow [\max(\mathcal{N}/c_0, 0)]^{2/3}. \tag{12}$$

In the classical almost sure Bahadur representation theory for sample quantiles, one has the error bound $O[n^{-3/4}(\log n)^{1/2}(\log \log n)^{1/4}]$ (see Shorack and Wellner [17]). We expect that the bound $O_{\mathbb{P}}(n^{-3/4})$ in (11) is optimal up to a multiplicative logarithmic factor.

Theorem 3(i) gives asymptotic properties of FDP, which is the value of the false discovery process $\Gamma_n$ at a random time $\nu_{\rm BH}$, while Theorem 3(ii) concerns false nondiscovery proportion (FNP). FNP is the proportion of null hypotheses being accepted which are actually false. Since $G$ is continuous, the FDR is $\alpha\pi_0$ (BH [2]). As pointed out in Genovese and Wasserman [9], it is not easy to study FDP since the random time $\nu_{\rm BH}$ and the false discovery process $\Gamma_n(\cdot)$ are dependent; recall (9). The relation (13) gives an asymptotic expansion for $\Gamma_n(\nu_{\rm BH}) - \alpha\pi_0$ with a good error bound $O_{\mathbb{P}}(n^{-3/4})$ and the term $\Lambda_n(\nu_0) - \Lambda(\nu_0)$ is easier to work with. It seems that the asymptotic expansion is new even in the special case of independent null hypotheses.

THEOREM 3. *Assume Condition 1 and $\alpha_*^{-1} > \alpha^{-1} > f(\nu_0)$.*

(i) *We have*

$$\Gamma_n(\nu_{\rm BH}) - \alpha\pi_0 = \frac{\alpha}{\nu_0}[\Lambda_n(\nu_0) - \Lambda(\nu_0)] + O_{\mathbb{P}}(n^{-3/4}). \tag{13}$$

*Consequently $\sqrt{n}[\Gamma_n(\nu_{\rm BH}) - \alpha\pi_0] \Rightarrow N(0, \sigma_0^2)$ for some $\sigma_0^2 < \infty$.*

(ii) *Let $X_n^* = \max_{i \leq n} X_i$ and define the false nondiscovery process*

$$\Xi_n(t) = \frac{\tilde{\Delta}_n(t)}{1 - F_n(t) + \mathbf{1}_{X_n^* \leq t}/n} \qquad \text{where } \tilde{\Delta}_n(t) = \frac{1}{n}\sum_{i=1}^n H_i \mathbf{1}_{X_i > t}.$$

*Let $c = \pi_0(\alpha - 1)/[1 - \alpha f(\nu_0)] + 1 - \nu_0/\alpha$ and $\Xi(t) = \pi_1[1 - G(t)]/[1 - F(t)]$. Then*

$$\Xi_n(\nu_{\rm BH}) - \Xi(\nu_0) = \frac{c[F_n(t) - F(t)]}{(1 - \nu_0/\alpha)^2} + \frac{\Delta_n(t) - \mathbb{E}\Delta_n(t)}{1 - \nu_0/\alpha} + O_{\mathbb{P}}(n^{-3/4}), \tag{14}$$

*and consequently $\sqrt{n}[\Xi_n(\nu_{\rm BH}) - \Xi(\nu_0)] \Rightarrow N(0, \sigma_1^2)$ for some $\sigma_1^2 < \infty$.*

We shall now discuss the estimation of the proportion of false null hypotheses $\pi_1$ and $g$ under dependence. When the $H_i$'s are independent, Genovese and Wasserman [9] pointed out that there is an unidentifiability issue in estimating $\pi_1$ and $g$ from the $p$-values $X_1, \ldots, X_n$. To see this, let $\lambda \in (1 - \min_{0 \leq x \leq 1} g(x), 1/\pi_1)$, $\pi_1^* = \lambda \pi$ and $g^*(x) = (g(x) - 1)/\lambda + 1$. Then we have the identity $(1 - \pi_1) + \pi_1 g(x) = (1 - \pi_1^*) + \pi_1^* g^*(x)$, suggesting that $X_i$ can also be viewed as a simple random sample from a mixture model with the two components: uniform$[0,1]$ and $g^*$. To ensure identifiability, we assume $g(1) = 0$. Since $f(t) = \pi_0 + \pi_1 g(t)$, as in Storey [18], we estimate $\pi_0 = 1 - \pi_1$ by

$$\hat{\pi}_0 = \frac{1 - F_n(1-b)}{b} \qquad \text{where } 0 < b < 1. \tag{15}$$

If $f$ is differentiable at 1, then in the sense of mean squared error the optimal bandwidth $b = b_n \asymp n^{-1/3}$ (cf. Lemma 2). Let $\hat{\pi}_1 = 1 - \hat{\pi}_0$ be the estimator of $\pi_1 = \mathbb{P}(H_j = 1)$. The BH procedure can be improved by the plug-in procedure: let

$$\nu_{\text{PI}} = \sup\{t \in [0,1] : t\hat{\pi}_0/\alpha \leq F_n(t)\}$$

and reject hypotheses for which $X_{(i)} \leq X_{(R_{\text{PI}})}$, where $R_{\text{PI}} = \lfloor n\nu_{\text{PI}}\hat{\pi}_0/\alpha \rfloor$. We argue that in the case of dependence null hypotheses, the plug-in procedure also improves the BH procedure by increasing power while it still controls the false discovery rate. Let

$$\nu_* = \sup\{t \in [0,1] : t\pi_0/\alpha \leq F(t)\}.$$

THEOREM 4. *Assume Condition 1, $g(1) = 0$ and $\alpha^{-1} > f(\nu_0)$. Further assume $\alpha/\pi_0 > \alpha_*$ and $b_n \asymp n^{-1/3}$. Then we have* (i)

$$\nu_{\text{PI}} - \nu_* = \frac{\nu_*(\pi_0 - \hat{\pi}_0)}{\pi_0 - \alpha f(\nu_*)} + \frac{F_n(\nu_*) - F(\nu_*)}{\pi_0/\alpha - f(\nu_*)} + O_{\mathbb{P}}(n^{-2/3}) \tag{16}$$

*and* (ii) $\Gamma_n(\nu_{\text{PI}}) - \alpha = \alpha(1 - \hat{\pi}_0/\pi_0) + O_{\mathbb{P}}(n^{-1/2})$.

**3. Examples and simulation studies.** Section 3.1 concerns one-dimensional processes and Section 3.2 contains an application to Ising models in $\mathbb{Z}^2$. In both cases we shall show that Condition 1 is satisfied.

3.1. *One-dimensional processes.* Assume that $(H_i)$ is a stationary process of the form

$$H_i = h(\ldots, \eta_{i-1}, \eta_i, \eta_{i+1}, \ldots), \tag{17}$$

where $\eta_i$ are i.i.d. random variables or innovations and $h$ is a measurable function. By allowing the dependence of $H_i$ on $\eta_j$, we are incorporating location information in modeling the dependence among null hypotheses. As a simple case, if $h$ in (17) is a function of $m$ ($m \in \mathbb{N}$) arguments: $H_i = h(\eta_{i-m+1}, \ldots, \eta_i)$, then $H_i$ is $m$-dependent. Our formulation (17) seems in line with the principle that "everything is related to everything else, but near things are more related than distant things" (Tobler [21]).

We now give a simple condition for the CLT $n^{-1/2}(N_n - n\pi_1) \Rightarrow N(0, \sigma^2)$, where $N_n = \sum_{i=1}^n H_i$. Let $\mathcal{F}_i = (\ldots, \eta_{i-1}, \eta_i)$ and define the projection operator $\mathcal{P}_k$ by $\mathcal{P}_k \xi = \mathbb{E}(\xi|\mathcal{F}_k) - \mathbb{E}(\xi|\mathcal{F}_{k-1})$ if the latter exists. Assume that

$$(18) \qquad c_0 := \sum_{i=-\infty}^{\infty} \delta_i < \infty \qquad \text{where } \delta_i = \|\mathcal{P}_0 H_i\|.$$

Then $\|N_n - n\pi_1\| \leq c_0 \sqrt{n}$ and the CLT holds (cf. Lemma 1). The quantity $\delta_i$ is related to the predictive dependence measure given in Wu [23]. Condition (18) indicates that the cumulative impact of $\eta_0$ in predicting the whole sequence $(H_i)_{i \in \mathbb{Z}}$ is finite. In this sense (18) is a short-range dependence condition (Wu [23]). If (18) is violated, then one enters the territory of long-range dependence and one may have a non-Gaussian limit.

We now verify (18) for truncation indicators of linear processes. Let $H_i = \mathbf{1}_{Z_i \leq z_*}$, where $z_* \in \mathbb{R}$ is fixed and $Z_k = \sum_{i=-\infty}^{\infty} a_i \eta_{k-i}$. Here $\eta_i$ are i.i.d. random variables and $(a_i)_{i \in \mathbb{Z}}$ are real coefficients. Let $f_\eta$ be the density of $\eta_i$ and $a_0 = 1$. Assume $\mathbb{E}(|\eta_i|^d) < \infty$, $d > 0$, and $c_* = \sup_z |f_\eta(z)| < \infty$. Let $d' = \min(1, d)$. Then $\delta_i = O(|a_i|^{d'/2})$ and (18) holds if $\sum_{i \in \mathbb{Z}} |a_i|^{d'/2} < \infty$. To this end, for $i \neq 0$ let $Y_i = Z_i - a_i \eta_0$. Since $Y_i - \eta_i$ and $\eta_i$ are independent, the density $f_{Y_i}$ of $Y_i$ satisfies $f_{Y_i}(y) = \mathbb{E} f_\eta[y - (Y_i - \eta_i)] \leq c_*$. Let $F_{Y_i}$ be the distribution function of $Y_i$. Then for $i \neq 0$,

$$(19) \quad \begin{aligned} \mathbb{E}|\mathbf{1}_{Z_i \leq z_*} - \mathbf{1}_{Y_i \leq z_*}| &\leq \mathbb{E}[\mathbb{E}(\mathbf{1}_{z_* - |a_i \eta_0| \leq Y_i \leq z_* + |a_i \eta_0|}|\eta_0)] \\ &= \mathbb{E}[F_{Y_i}(z_* + |a_i \eta_0|) - F_{Y_i}(z_* - |a_i \eta_0|)] \\ &\leq \mathbb{E}\{\min(1, 2c_*|a_i \eta_0|)\} \\ &\leq \mathbb{E}\{(2c_*|a_i \eta_0|)^{d'}\} = O(|a_i|^{d'}). \end{aligned}$$

Let $\eta_0', \eta_i, i \in \mathbb{Z}$, be i.i.d. and $Z_i' = Y_i + a_i \eta_0'$. Then (19) implies $\mathbb{E}|\mathbf{1}_{Z_i \leq z_*} - \mathbf{1}_{Z_i' \leq z_*}| = O(|a_i|^{d'})$. Observe that $\mathbb{E}(\mathbf{1}_{Z_i \leq z_*} - \mathbf{1}_{Z_i' \leq z_*}|\mathcal{F}_0) = \mathcal{P}_0 H_i$. By Jensen's inequality, $\delta_i = O(|a_i|^{d'/2})$.

3.2. *Ising models.* Markov random fields have been widely used in image analysis and spatial statistics. Here we shall consider a false discovery control paradigm with the null hypotheses $(H_s)$ satisfying the Gibbs distribution in $\mathbb{Z}^2$ and thus $(H_s)$ are spatially dependent. Let $L_s = 2H_s - 1$. Then $L_s = -1$ (resp. 1) implies that the null hypothesis $H_s$ is true (resp. false). That $L_s = 1$ may imply that a neuron is excited or a plant is infected. Here we consider the simplest Ising model. For a site $s = (j, k) \in \mathbb{Z}^2$, let $\mathcal{N}_s = \{(j', k') \in \mathbb{Z}^2 : |j - j'| + |k - k'| = 1\}$ be the neighborhood of $s$ and write $\mathbb{Z}^2 \setminus s = \{t \in \mathbb{Z}^2 : t \neq s\}$. For a set $A \subset \mathbb{Z}^2$ write $L_A = (L_a, a \in A)$ and $l_A = (l_a, a \in A)$, where $l_a \in \{-1, 1\}$, $a \in \mathbb{Z}^2$. Assume that we have the Markovian structure

$$
\begin{aligned}
\mathbb{P}[L_s = l_s | L_{\mathbb{Z}^2 \setminus s} = l_{\mathbb{Z}^2 \setminus s}] &= \mathbb{P}[L_s = l_s | L_{\mathcal{N}_s} = l_{\mathcal{N}_s}] \\
&= \frac{\exp(\beta l_s \sum_{t \in \mathcal{N}_s} l_t)}{\exp(\beta \sum_{t \in \mathcal{N}_s} l_t) + \exp(-\beta \sum_{t \in \mathcal{N}_s} l_t)},
\end{aligned}
\tag{20}
$$

where $\beta$ characterizes the interaction between pairs of nearest-neighbor spins and it is a function of Boltzmann's constant and the temperature. That $\beta > 0$ (resp. $\beta < 0$) corresponds to ferromagnetic (resp. antiferromagnetic) interaction. The former is an attractive feature in dealing with situations in which one expects that false null hypotheses occur in clumps or clusters. In the antiferromagnetic case, one has negative dependence which prevents false null hypotheses from occurring in clumps.

With (20), the distribution of $H_s$ only depends on the values of $H$ at the four neighbors of $s$. For more details see Winkler [22]. Let $\beta_* = 2^{-1} \log(1 + \sqrt{2}) = 0.4406868\ldots$ be the critical value. If $0 \leq \beta < \beta_*$, then $\mathbb{E}(L_s) = 0$, and we can apply the central limit theorem in Newman [14] or Baker and Krinsky [1]: the covariance $\text{cov}(L_0, L_s) \to 0$ decays to zero exponentially quickly as $|s| \to \infty$ and $n^{-1/2}(N_n - n\pi_1) \Rightarrow N(0, \sigma^2)$. So Condition 1 is satisfied and Theorems 1–4 are applicable.

Consider the situation that $(H_s)$ are not directly observable and we want to test whether $H_s = 0$ or $H_s = 1$. We conduct pixel-wise multiple hypothesis tests. Assume that for each site $s$, under $H_s = 0$, the $p$-value $X_s$ has a uniform$(0, 1)$ distribution while $[X_s | H_s = 1] \sim G$. Since the underlying $(H_i)$ is not observed and one only knows $p$-values $X_i$ which are calculated from test statistics, we are thus dealing with hidden Markov models by viewing $(H_i)$ as hidden states. Analysis of the $p$-value sequence $(X_i)$ is useful in understanding the dependence structure of $(H_i)$ and provides spatial information of false null hypotheses.

In our simulation we choose the lattice set $\{1, 2, \ldots, 50\}^2$ with periodic boundary conditions and choose seven levels of $\beta$: $\beta = -0.3$, 0, 0.1, 0.2, 0.3, 0.4 and 0.44. Note that $\beta = 0$ implies independent null hypotheses and larger $\beta$ indicates stronger dependence. The density of the alternative distribution is $g(x) = a(1+a)^2/(x+a)^2 - a$, $x \in (0, 1)$, where $a = 1/98$. Then $g(1) = 0$, $g(0) = 100$ and the quantity $\alpha_*$ in (10) is $2/101$.

Our simulation study shows that, if the dependence is relatively weaker, then $\Gamma_n(\nu_{\text{BH}})$ is more concentrated on $\alpha \pi_0$ and the approximation (13) in Theorem 3 is better. We apply the Gibbs sampler with random sweeps (Greenwood, McKeague and Wefelmeyer [10]) and the number of iterations is $1.25 \times 10^6$. Choose the level $\alpha = 0.1$. For $\beta < \beta_*$, we have $\pi_0 = 1/2$ and $\alpha\pi_0 = 0.05$. Write $\delta_1 = \Gamma_n(\nu_{\mathrm{BH}}) - \alpha\pi_0$ and $\delta_2 = \nu_0^{-1}\alpha[\Lambda_n(\nu_0) - \Lambda(\nu_0)]$. Table 1 shows the estimated $\mathbb{E}(\delta_1^2)$ and $\mathbb{E}(|\delta_1 - \delta_2|^2)$ based on 100 repetitions. It suggests that $\delta_2$ approximates $\delta_1$ reasonably well. As the dependence gets stronger, $\mathbb{E}(\delta_1^2)$ becomes larger and the false discovery proportion $\Gamma_n(\nu_{\mathrm{BH}})$ is less concentrated on $\alpha\pi_0$.

Genovese, Lazar and Nichols [8] showed that the false discovery rate controlling procedure can be useful in the analysis of image data. Figure 1 shows image restoration based on the $p$-values under the conditional independence model. In our simulation we applied pixel-wise multiple hypothesis tests with FDR-controlling procedure and the level is $\alpha = 0.1$. The first row is the simulated Ising images for $\beta = 0.3$ and 0.44, respectively. The second row shows the estimated images and the third row gives the differences. The red (resp. blue) dots are false positives (resp. negatives). With larger $\alpha$ (say $\alpha = 0.15$), the number of false negatives is reduced (the simulation is not reported in the paper).

Figure 1 suggests that, if the dependence is strong (e.g., $\beta = 0.44$) and the false null hypotheses are clustered, then it is possible to improve the restored images by incorporating the spatial dependence structure. Pacifico et al. [15] applied FDR-thresholding to construct conservative confidence envelopes for Gaussian random fields.

**4. Proofs.** This section provides proofs of results stated in Section 2. For readability we list necessary notation here. Recall (6) and (7) for $\Lambda_n(t)$, $F_n(t)$ and $\Delta_n(t)$. Let $N_n = \sum_{i=1}^n H_i$ be the total number of false null hypotheses,

$$(21) \quad \Lambda_n^*(t) = \frac{1}{n}\sum_{i=1}^n (1 - H_i)t = t(1 - N_n/n) \quad \text{and} \quad \Delta_n^*(t) = G(t)N_n/n.$$

TABLE 1
*The estimated $\mathbb{E}(\delta_1^2)$ and $\mathbb{E}(|\delta_1 - \delta_2|^2)$ based on 100 repetitions*

| $\beta$ | $\hat{\mathbb{E}}(\delta_1^2)$ | $\hat{\mathbb{E}}(|\delta_1 - \delta_2|^2)$ |
|---|---|---|
| $-0.3$ | $4.2 \times 10^{-5}$ | $6.2 \times 10^{-7}$ |
| $0$ | $4.4 \times 10^{-5}$ | $1.1 \times 10^{-6}$ |
| $0.1$ | $5.6 \times 10^{-5}$ | $1.7 \times 10^{-6}$ |
| $0.2$ | $5.7 \times 10^{-5}$ | $1.5 \times 10^{-6}$ |
| $0.3$ | $6.0 \times 10^{-5}$ | $3.4 \times 10^{-6}$ |
| $0.4$ | $9.5 \times 10^{-5}$ | $7.6 \times 10^{-6}$ |
| $0.44$ | $7.1 \times 10^{-4}$ | $1.1 \times 10^{-4}$ |

Here $\delta_1 = \Gamma_n(\nu_{\mathrm{BH}}) - \alpha\pi_0$ and $\delta_2 = \nu_0^{-1}\alpha[\Lambda_n(\nu_0) - \Lambda(\nu_0)]$.

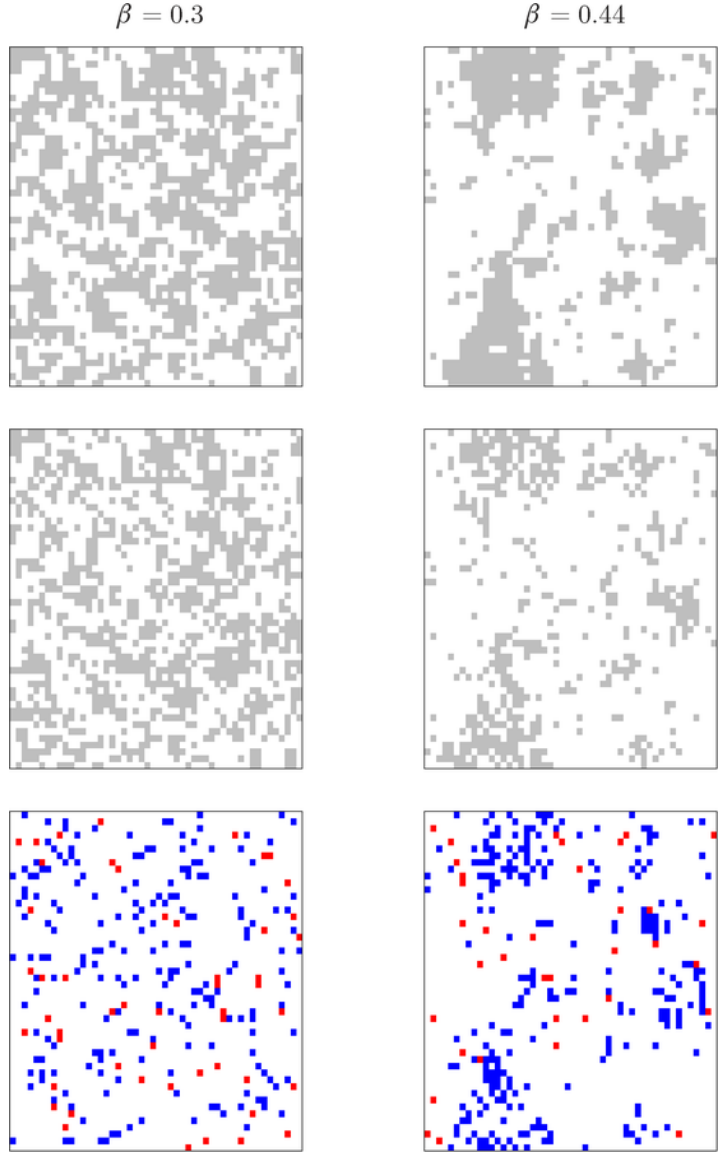

FIG. 1. *Row* 1: *simulated Ising models for $\beta = 0.3$ and $0.44$, respectively. Row* 2: *restored images based on the p-values under the conditional independence model. Here we applied pixel-wise multiple hypothesis tests with FDR-controlling procedure and the level is $\alpha = 0.1$. Row* 3: *the differences between the restored images and the original ones. Dots in red (resp. blue) are false positives (resp. negatives).*

Write $F_n^* = \Lambda_n^* + \Delta_n^*$. Define

$$\Omega_n(t) = \frac{1}{n} \sum_{i=1}^{n} (1 - H_i)(\mathbf{1}_{U_i \le t} - t), \tag{22}$$

$$\Xi_n(t) = \frac{1}{n} \sum_{i=1}^{n} H_i[\mathbf{1}_{U_i \leq G(t)} - G(t)]. \tag{23}$$

LEMMA 1. *Assume* (18). *Then* $\|N_n - n\pi_1\| \leq c_0\sqrt{n}$ *and* $n^{-1/2}(N_n - n\pi_1) \Rightarrow N(0, \sigma^2)$.

PROOF. By stationarity, $\|\mathcal{P}_k N_n\| \leq \sum_{i=1}^{n} \delta_{i-k} \leq c_0$. Since $\mathcal{P}_k$ are orthogonal, we have

$$\|N_n - n\pi_1\|^2 = \sum_{k \in \mathbb{Z}} \|\mathcal{P}_k N_n\|^2 \leq \sum_{k \in \mathbb{Z}} c_0 \sum_{i=1}^{n} \delta_{i-k} = nc_0^2.$$

A similar version of the CLT is given in Hannan [11] and the argument therein is applicable here. Let $D_k = \sum_{i \in \mathbb{Z}} \mathcal{P}_k H_i$ and $M_n = \sum_{k=1}^{n} D_k$. Then $D_k$ are stationary martingale differences. Let $u_j = \sum_{i=j}^{\infty} \delta_i$ and $l_j = \sum_{i=-\infty}^{j} \delta_i$. Since $\mathcal{P}_k$, $k \in \mathbb{Z}$, are orthogonal,

$$\|N_n - n\pi_1 - M_n\|^2 = \sum_{k \in \mathbb{Z}} \|\mathcal{P}_k(N_n - M_n)\|^2.$$

If $k \leq 0$, then $\|\mathcal{P}_k(N_n - M_n)\| = \|\mathcal{P}_k N_n\| \leq \sum_{i=1}^{n} \delta_{i-k}$. So $\sum_{k=-\infty}^{0} \|\mathcal{P}_k N_n\|^2 \leq c_0 \sum_{i=1}^{n} u_i = o(n)$ since $u_m \to 0$ as $m \to \infty$. Similarly, $\sum_{k=n+1}^{\infty} \|\mathcal{P}_k N_n\|^2 = o(n)$. For $1 \leq k \leq n$, since $\mathcal{P}_k M_n = D_k$, $\|\mathcal{P}_k(N_n - M_n)\| \leq u_{n+1-k} + l_{-k}$. So we also have $\sum_{k=1}^{n} \|\mathcal{P}_k(N_n - M_n)\|^2 = o(n)$. Thus $\|N_n - n\pi_1 - M_n\|^2 = o(n)$. By the martingale CLT, $M_n/\sqrt{n} \Rightarrow N(0, \sigma^2)$ with $\sigma = \|D_k\|$. So the lemma holds. □

LEMMA 2. *Assume* $\sup_{x \in [0,1]} g(x) < \infty$. *Let $b_n$ be a sequence of bandwidths satisfying*

$$b_n \to 0 \quad \text{and} \quad nb_n \to \infty. \tag{24}$$

*Then under Condition* 1, *we have*

$$\sqrt{n/b_n}[F_n(b_n) - F(b_n)] \Rightarrow N(0, f(0)) \tag{25}$$

*and*

$$\sqrt{n/b_n}[F_n(1-b_n) - F(1-b_n)] \Rightarrow N(0, f(1)). \tag{26}$$

PROOF. Denote by $\sqrt{-1}$ the imaginary unit. Let

$$D_i = (1 - H_i)(\mathbf{1}_{U_i \leq b_n} - b_n) + H_i(\mathbf{1}_{U_i \leq G(b_n)} - G(b_n))$$

and $Q_n = \sum_{i=1}^n D_i$. Let $t \in \mathbb{R}$ be a fixed number and $t_n = t/\sqrt{nb_n}$. Under the conditional independence model (4), the conditional characteristic function

$$\begin{aligned}\phi_n(t) &:= \mathbb{E}[\exp(\sqrt{-1}t_n Q_n)|H_i, 1 \leq i \leq n]\\ &= [G(b_n)\exp(\sqrt{-1}t_n(1-G(b_n)))\\ &\quad + (1-G(b_n))\exp(-\sqrt{-1}t_n G(b_n))]^{N_n}\\ &\quad + [b_n \exp(\sqrt{-1}t_n(1-b_n)) \times (1-b_n)\exp(-\sqrt{-1}t_n b_n)]^{n-N_n}.\end{aligned}$$

By Condition 1, $N_n/n \to \pi_1$ in probability. Using Taylor's expansions $\exp(\delta) = 1 + \delta + \delta^2/2 + O(\delta^3)$, $G(\delta) = \delta g(0) + o(\delta)$, after elementary calculations we have

$$\phi_n(t) \to \exp\{-t^2/[2\pi_0 + 2\pi_1 g(0)]\} = \exp\{-t^2/[2f(0)]\} \qquad \text{in probability.}$$

By the Lebesgue dominated convergence theorem, $\mathbb{E}[\phi_n(t)] \to \exp\{-t^2/[2f(0)]\}$ since $|\phi_n(t)| \leq 1$. So $Q_n/\sqrt{nb_n} \Rightarrow N(0, f(0))$. By Condition 1, since $G(b_n) = O(b_n)$, we have

$$|n[F_n(b_n) - F(b_n)] - Q_n| \leq [b_n + G(b_n)]|N_n - n\pi_1| = O_{\mathbb{P}}(b_n\sqrt{n}).$$

So (25) follows since $b_n \to 0$. The other assertion (26) can be similarly proved by considering $D_i' = (1-H_i)(\mathbf{1}_{U_i > 1-b_n} - b_n) + H_i(\mathbf{1}_{U_i > 1-G(b_n)} - G(b_n))$. □

LEMMA 3. *Let $b_n$ be a sequence of positive numbers satisfying $b_n \in (0,1)$ and $nb_n \to \infty$. Assume that $\sup_{x \in [0,1]} g(x) < \infty$. Then we have*

$$\sup_{|u| \leq b_n} n|\Omega_n(t+u) - \Omega_n(t)| = O_{\mathbb{P}}[(nb_n)^{1/2}] \tag{27}$$

*and*

$$\sup_{|u| \leq b_n} n|\Xi_n(t+u) - \Xi_n(t)| = O_{\mathbb{P}}[(nb_n)^{1/2}]. \tag{28}$$

PROOF. For i.i.d. uniform$(0,1)$ random variables $U_i$, $i \in \mathbb{Z}$, let $W_n(u) = \sum_{i=1}^n \mathbf{1}_{U_i \leq u} - nu$. By Lemma 2.3 in Stute [20], there exists a constant $c_0$ such that

$$\mathbb{P}\left[\sup_{0 \leq u \leq b} |W_n(u+t) - W_n(t)| > s\sqrt{nb}\right] \leq 4e^{-s^2/16} \tag{29}$$

holds for all $0 < b < 1/8$ and $32 \leq s \leq c_0\sqrt{nb}$. Since $(H_i)$ is 0/1-valued and it is independent of $U_i$, it is easily seen that (29) implies

$$\mathbb{P}\left[\sup_{0 \leq u \leq b} |n\Omega_n(t+u) - n\Omega_n(t)| > s\sqrt{nb}\right] \leq 4e^{-s^2/16}.$$

So we have (27) since $nb_n \to \infty$. A similar argument entails (28). □

LEMMA 4. *Assume Condition 1 and $\sup_{x\in[0,1]} g(x) < \infty$. Let $b_n \in (0,1)$. Then for $\Lambda_n^*(t)$ and $\Delta_n^*(t)$ defined in Lemma 3, we have*

$$(30) \quad \sup_{|u|\leq b_n} n|[\Lambda_n^*(t+u) - \Lambda(t+u)] - [\Lambda_n^*(t) - \Lambda(t)]| = O_{\mathbb{P}}(b_n n^{1/2})$$

*and*

$$(31) \quad \sup_{|u|\leq b_n} n|[\Delta_n^*(t+u) - \Delta(t+u)] - [\Delta_n^*(t) - \Delta(t)]| = O_{\mathbb{P}}(b_n n^{1/2}).$$

PROOF. Since $\Delta_n^*(t) = G(t)N_n/n$ and $\sup_{x\in[0,1]} g(x) < \infty$, by Condition 1, we have (31). Similarly (30) follows. □

4.1. *Proof of Theorem 1.* By the weak convergence theory (Billingsley [4]), it suffices to establish (i) the finite-dimensional convergence and (ii) the tightness.

We first show that the process $\sqrt{n}\{\Lambda_n(t) - \Lambda(t)\}$ is tight. Let $L_n(t) = \Lambda_n^*(t) - \Lambda(t) = t(\pi_1 - N_n/n)$, $0 \leq t \leq 1$. By Condition 1, $\sqrt{n}(\pi_1 - N_n/n) = O_{\mathbb{P}}(1)$. So $\sqrt{n}L_n(t)$ is trivially tight. Note that $\Lambda_n(t) - \Lambda(t) = \Omega_n(t) + L_n(t)$. Following the tightness argument for the process $n^{-1/2}\sum_{i=1}^n (\mathbf{1}_{U_i\leq u} - u)$, $0 \leq u \leq 1$ (cf. Theorem 16.4 in Billingsley [4]), since $(H_i)$ and $(U_i)$ are independent, we can easily derive that the process $\sqrt{n}\Omega_n(t)$ is also tight. Similarly, we can show that $\sqrt{n}\{\Lambda_n(t) - \Lambda(t)\}$ is tight by noting that $\sup_{0\leq t\leq 1} g(t) < \infty$. So $(\sqrt{n}\{\Lambda_n(t) - \Lambda(t)\}, \sqrt{n}\{\Delta_n(t) - \Delta(t)\})$ is tight.

We now show the finite-dimensional convergence. Let $a, b$ be two real numbers; let

$$(32) \quad T_n = \sum_{i=1}^n (J_i - \mathbb{E}J_i) \quad \text{where } J_i = a(1-H_i)\mathbf{1}_{U_i\leq t} + bH_i\mathbf{1}_{U_i\leq G(t)}.$$

We shall calculate the characteristic function $\varphi_n(\theta) = \mathbb{E}\{\exp[\theta\sqrt{-1}T_n/\sqrt{n}]\}$, $\theta \in \mathbb{R}$. Let $A(\theta) = \log\mathbb{E}\{\sqrt{-1}\theta\mathbf{1}_{U_1\leq t}\}$ and $B(\theta) = \log\mathbb{E}\{\sqrt{-1}\theta\mathbf{1}_{G(U_1)\leq t}\}$. Then for small $|\delta|$, we have

$$A(\delta) = \log(1 - t + te^{\sqrt{-1}\delta}) = t\delta\sqrt{-1} - \frac{\delta^2}{2}t(1-t) + O(\delta^3),$$

$$B(\delta) = \log\{1 - G(t) + G(t)e^{\sqrt{-1}\delta}\}$$
$$= G(t)\delta\sqrt{-1} - \frac{\delta^2}{2}G(t)[1 - G(t)] + O(\delta^3).$$

Let $v = G(t)\theta b - t\theta a$, $\varrho_0 = t(1-t)\theta^2 a^2/2$ and $\varrho_1 = G(t)[1-G(t)]\theta^2 b^2/2$. With the preceding two relations, since $(N_n - n\pi_1)/\sqrt{n} \Rightarrow N(0, \sigma^2)$, as the argument for $\phi_n(t)$ in the proof of Lemma 2, we have

$$
\begin{aligned}
\lim_{n \to \infty} \varphi_n(\theta) &= \lim_{n \to \infty} \frac{\mathbb{E} \exp\{(n - N_n)A(\theta a/\sqrt{n}) + N_n B(\theta b/\sqrt{n})\}}{\exp\{\sqrt{-1}\theta\sqrt{n}[a\pi_0 t + b\pi_1 G(t)]\}} \\
&= \lim_{n \to \infty} \mathbb{E} \exp\{(N_n - n\pi_1)\sqrt{-1}v/\sqrt{n} \\
&\quad - (1 - N_n/n)\varrho_0 - (N_n/n)\varrho_1\} \\
&= \mathbb{E} \exp\{-v^2\sigma^2/2 - \pi_0\varrho_0 - \pi_1\varrho_1\}
\end{aligned}
\tag{33}
$$

after elementary manipulations. Hence $T_n/\sqrt{n}$ is asymptotically normal. Consequently, by the Crámer–Wold device, the finite-dimensional convergence follows.

4.2. *Proof of Theorem 2.* (i) Let $b_n$ be a real sequence with $b_n \in (0,1)$ and $nb_n \to \infty$. Since $b_n \leq \sqrt{b_n}$, by Lemmas 3 and 4, we have

$$
\sup_{|u| \leq b_n} n|[F_n(t+u) - F(t+u)] - [F_n(t) - F(t)]| = O_\mathbb{P}((nb_n)^{1/2}). \tag{34}
$$

We first show that $\sqrt{n}(\nu_{\text{BH}} - \nu_0) = O_\mathbb{P}(1)$. To this end, it suffices to show that for any positive sequence $B_n \to \infty$, $\sqrt{n}(\nu_{\text{BH}} - \nu_0) = O_\mathbb{P}(B_n)$. Without loss of generality assume $B_n \leq \log n$ since otherwise we can let $B_n' = \min(B_n, \log n)$. Applying (34) with $b_n = B_n/\sqrt{n}$, since $F_n(\nu_0) - F(\nu_0) = O_\mathbb{P}(n^{-1/2})$, we have

$$
\begin{aligned}
F_n(\nu_0 + b_n) &= F(\nu_0 + b_n) + [F_n(\nu_0) - F(\nu_0)] + O_\mathbb{P}((b_n/n)^{1/2}) \\
&= F(\nu_0) + b_n f(\nu_0) + O(b_n^2) + O_\mathbb{P}(n^{-1/2}).
\end{aligned}
\tag{35}
$$

Note that $t > \nu_{\text{BH}}$ if and only if $t/\alpha > F_n(t)$. Since $1/\alpha > f(\nu_0)$, $F(\nu_0) = \nu_0/\alpha$ and $B_n \to \infty$, we have by (35) that

$$\mathbb{P}(\nu_0 + b_n > \nu_{\text{BH}}) = \mathbb{P}[(\nu_0 + b_n)/\alpha > F_n(\nu_0 + b_n)] \to 1$$

as $n \to \infty$. Similarly, $\mathbb{P}(\nu_0 - b_n \leq \nu_{\text{BH}}) \to 1$. So $\sqrt{n}(\nu_{\text{BH}} - \nu_0) = O_\mathbb{P}(1)$, which, by another application of (34) with $b_n = C/\sqrt{n}$, implies

$$
n|[F_n(\nu_{\text{BH}}) - F(\nu_{\text{BH}})] - [F_n(\nu_0) - F(\nu_0)]| = O_\mathbb{P}(n^{1/4}). \tag{36}
$$

Since $|F_n(\nu_{\text{BH}}) - \nu_{\text{BH}}/\alpha| = O(n^{-1})$ and $F(\nu_{\text{BH}}) = F(\nu_0) + (\nu_{\text{BH}} - \nu_0)f(\nu_0) + O_\mathbb{P}(n^{-1})$, (11) follows.

The CLT $\sqrt{n}(\nu_{\text{BH}} - \nu_0) \Rightarrow N(0, \sigma^2)$ easily follows from (11) in view of the argument of (32) and (33) in the proof of Theorem 1: let $a = b = 1$ in (32), then $J_i = \mathbf{1}_{X_i \leq t}$.

(ii) As in (i) we shall show that for any positive sequence $B_n \to \infty$, $R = O_\mathbb{P}(B_n)$. To this end, let $b_n = B_n/n$, $t_n = n[b_n - F(\alpha b_n)]/\sqrt{nb_n}$ and

$$Z_n = \frac{n(F_n(\alpha b_n) - F(\alpha b_n))}{\sqrt{nb_n}}.$$

Since $b_n \to 0$ and $nb_n = B_n \to \infty$, by Taylor's expansion, $F(\alpha b_n) = f(0)\alpha b_n + o(b_n)$. Hence $t_n/\sqrt{nb_n} \to 1 - \alpha/\alpha_* > 0$. So $t_n \to \infty$. By Lemma 2, $Z_n \Rightarrow N[0, \alpha f(0)]$. Therefore
$$\mathbb{P}(R < B_n) = \mathbb{P}[F_n(\alpha b_n) < b_n] = \mathbb{P}(t_n > Z_n) \to 1.$$

(iii) Let $z > 0$ be fixed and $b_n = n^{-1/3}z$. By Taylor's expansion, $F(b_n) = b_n f(0) + b_n^2 f'(0)/2 + o(b_n^2)$. Hence $u_n := \sqrt{n/b_n}[b_n/\alpha - F(b_n)] \to -f'(0)z^{3/2}/2$. By Lemma 2(i),
$$\begin{aligned}\mathbb{P}(\nu_{\text{BH}} < b_n) &= \mathbb{P}[F_n(b_n) < b_n/\alpha] \\ &= \mathbb{P}\{\sqrt{n/b_n}[F_n(b_n) - F(b_n)] < u_n\} \\ &\to \mathbb{P}\left\{\sqrt{f(0)}\mathcal{N} \leq -\frac{f'(0)}{2}z^{3/2}\right\} \\ &= \mathbb{P}\{[\max(\mathcal{N}/c_0, 0)]^{2/3} \leq z\},\end{aligned}$$
which proves (12).

4.3. *Proof of Theorem* 3. (i) By Theorem 2(i), $\nu_{\text{BH}} - \nu_0 = O_\mathbb{P}(1/\sqrt{n})$. Similarly as in the proof of (36), by Lemmas 3 and 4, we have

(37) $\quad \Lambda_n(\nu_{\text{BH}}) = \Lambda_n(\nu_0) + \Lambda(\nu_{\text{BH}}) - \Lambda(\nu_0) + O_\mathbb{P}(n^{-3/4}).$

Recall $\Lambda(t) = t\pi_0$. Observe that $F(\nu_{\text{BH}}) - F(\nu_0) = (\nu_{\text{BH}} - \nu_0)f(\nu_0) + O_\mathbb{P}(1/n)$, $F_n(\nu_0) = F(\nu_0) + O_\mathbb{P}(1/\sqrt{n})$ and, by (36), $F_n(\nu_{\text{BH}}) = F(\nu_0) + O_\mathbb{P}(1/\sqrt{n})$. By (36) and (37), we have

(38) $\quad \Lambda_n(\nu_{\text{BH}})F(\nu_0) - F_n(\nu_{\text{BH}})\Lambda(\nu_0) = \frac{1}{n}\sum_{i=1}^n \{J_i - \mathbb{E}(J_i)\} + O_\mathbb{P}(n^{-3/4}),$

where
$$\begin{aligned}J_i &= F(\nu_0)(1 - H_i)\mathbf{1}_{X_i \leq \nu_0} - \Lambda(\nu_0)\mathbf{1}_{X_i \leq \nu_0} \\ &\quad + \mathbf{1}_{X_i \leq \nu_0}\frac{F(\nu_0)\Lambda'(\nu_0) - f(\nu_0)\Lambda(\nu_0)}{\alpha^{-1} - f(\nu_0)}.\end{aligned}$$

Since $F(\nu_0) = \nu_0/\alpha$ and $\Lambda'(\nu_0) = \pi_0$, simple calculations show $J_i = F(\nu_0)(1 - H_i)\mathbf{1}_{X_i \leq \nu_0}$. Note that $F(\nu_0) < \nu_0$ and $F(\nu_0) = \pi_0\nu_0 + \pi_1 G(\nu_0)$. Then $G(\nu_0) > \nu_0$. Using the property of conditional independence,
$$\mathbb{P}\left(\min_{i \leq n} X_i \geq \nu_0 | H_1, \ldots, H_n\right) = (1 - \nu_0)^{n-N_n}(1 - G(\nu_0))^{N_n} \leq (1 - \nu_0)^n.$$
So $\mathbb{P}(\min_{i \leq n} X_i \geq \nu_0) \leq (1 - \nu_0)^n$ and hence (13) follows from (38) by noting that $F_n(\nu_{\text{BH}}) = F(\nu_0) + O_\mathbb{P}(1/\sqrt{n})$. The CLT $\sqrt{n}[\Gamma_n(\nu_{\text{BH}}) - \alpha\pi_0] \Rightarrow N(0, \sigma_0^2)$ follows from (33).

(ii) The argument is similar to the one in (i). We have an analog of (37) with $\Lambda_n(\cdot)$ therein replaced by $\tilde{\Lambda}_n(\cdot)$ and (14) similarly holds. The CLT also follows from (33). Details are omitted.

4.4. *Proof of Theorem* 4. For (16), the argument in the proof of Theorem 3 is applicable. Let $B_n$ be a positive sequence that diverges to infinity slower than $\log n$; let $r_n = b_n + (nb_n)^{-1/2} \asymp n^{-1/3}$. By Lemmas 3 and 4,

$$
\begin{aligned}
F_n(\nu_* + r_n B_n) &= F(\nu_* + r_n B_n) \\
&\quad + [F_n(\nu_*) - F(\nu_*)] + O_{\mathbb{P}}[(r_n B_n/n)^{1/2}] \\
&= F(\nu_*) + f(\nu_*) r_n B_n + O(r_n^2 B_n^2) + O_{\mathbb{P}}(n^{-1/2}).
\end{aligned}
\tag{39}
$$

By Lemma 2, $\sqrt{nb_n}(\hat{\pi}_0 - \mathbb{E}\hat{\pi}_0) \Rightarrow N(0, f(1))$. Since $b_n \asymp n^{-1/3}$ and $B_n \to \infty$, we have

$$(40) \quad \mathbb{P}\{(\nu_* + r_n B_n)(\hat{\pi}_0 - \mathbb{E}\hat{\pi}_0) \geq \nu_*(\pi_0 - \mathbb{E}\hat{\pi}_0) + r_n B_n[\alpha f(\nu_*) - \mathbb{E}\hat{\pi}_0]\} \to 1$$

since $f(\nu_*) < \pi_0/\alpha$ and $\mathbb{E}\hat{\pi}_0 = \pi_0 + O(b_n)$. Note that $F(\nu_*) = \pi_0 \nu_*/\alpha$. By (39) and (40),

$$\mathbb{P}[(\nu_* + r_n B_n)\hat{\pi}_0 > \alpha F_n(\nu_* + r_n B_n)] \to 1,$$

which implies that $\mathbb{P}(\nu_{\mathrm{PI}} \leq \nu_* + r_n B_n) \to 1$. Similarly, we have $\mathbb{P}(\nu_{\mathrm{PI}} \geq \nu_* - r_n B_n) \to 1$ and hence $\nu_{\mathrm{PI}} - \nu_* = O_{\mathbb{P}}(r_n)$. By Lemmas 3 and 4,

$$(41) \quad F_n(\nu_{\mathrm{PI}}) = F(\nu_{\mathrm{PI}}) + [F_n(\nu_*) - F(\nu_*)] + O_{\mathbb{P}}[(r_n/n)^{1/2}].$$

Since $|F_n(\nu_{\mathrm{PI}}) - \nu_{\mathrm{PI}}\hat{\pi}_0/\alpha| \leq n^{-1}$ and $F(\nu_{\mathrm{PI}}) - F(\nu_*) = (\nu_{\mathrm{PI}} - \nu_*)f(\nu_*) + O_{\mathbb{P}}(r_n^2)$, (16) follows from (41) after elementary calculations.

Using the argument in the proof of Theorem 3, we can similarly obtain (ii) with no essential difficulties. Since the calculation is lengthy, the details are omitted.

**Acknowledgments.** We would like to thank the referee and an Associate Editor for their many helpful comments.

## REFERENCES

[1] BAKER, JR., G. A. and KRINSKY, S. (1977). Renormalization group structure for translationally invariant ferromagnets. *J. Math. Phys.* **18** 590–607. MR0443771
[2] BENJAMINI, Y. and HOCHBERG, Y. (1995). Controlling the false discovery rate: A practical and powerful approach to multiple testing. *J. Roy. Statist. Soc. Ser. B* **57** 289–300. MR1325392
[3] BENJAMINI, Y. and YEKUTIELI, D. (2001). The control of the false discovery rate in multiple testing under dependency. *Ann. Statist.* **29** 1165–1188. MR1869245
[4] BILLINGSLEY, P. (1968). *Convergence of Probability Measures*. Wiley, New York. MR0233396
[5] CHI, Z. (2005). Criticality of a false discovery rate controlling procedure. Technical Report 05-25, Dept. Statistics, Univ. Connecticut.
[6] DE CASTRO, M. C. and SINGER, B. H. (2006). Controlling the false discovery Rate: A new application to account for multiple and dependent tests in local statistics of spatial association. *Geographical Analysis* **38** 180–208.

[7] FARCOMENI, A. (2004). Multiple testing procedures under dependence, with applications. Ph.D. thesis, Univ. Roma "La Sapienza."
[8] GENOVESE, C., LAZAR, N. A. and NICHOLS, T. (2002). Thresholding of statistical maps in functional neuroimaging using the false discovery rate. *Neuroimage* **15** 870–878.
[9] GENOVESE, C. and WASSERMAN, L. (2004). A stochastic process approach to false discovery control. *Ann. Statist.* **32** 1035–1061. MR2065197
[10] GREENWOOD, P. E., MCKEAGUE, I. W. and WEFELMEYER, W. (1998). Information bounds for Gibbs samplers. *Ann. Statist.* **26** 2128–2156. MR1700224
[11] HANNAN, E. J. (1973). Central limit theorems for time series regression. *Z. Wahrsch. Verw. Gebiete* **26** 157–170. MR0331683
[12] MCLACHLAN, G. and PEEL, D. (2000). *Finite Mixture Models.* Wiley, New York. MR1789474
[13] MEINSHAUSEN, N. and RICE, J. (2006). Estimating the proportion of false null hypotheses among a large number of independently tested hypotheses. *Ann. Statist.* **34** 373–393. MR2275246
[14] NEWMAN, C. M. (1980). Normal fluctuations and the FKG inequalities. *Comm. Math. Phys.* **74** 119–128. MR0576267
[15] PACIFICO, M. P., GENOVESE, C., VERDINELLI, I. and WASSERMAN, L. (2004). False discovery control for random fields. *J. Amer. Statist. Assoc.* **99** 1002–1014. MR2109490
[16] QIU, X., BROOKS, A. I., KLEBANOV, L. and YAKOVLEV, N. (2005). The effects of normalization on the correlation structure of microarray data. *BMC Bioinformatics* **6** 120.
[17] SHORACK, G. R. and WELLNER, J. A. (1986). *Empirical Processes with Applications to Statistics.* Wiley, New York. MR0838963
[18] STOREY, J. D. (2002). A direct approach to false discovery rates. *J. R. Stat. Soc. Ser. B Stat. Methodol.* **64** 479–498. MR1924302
[19] STOREY, J. D., TAYLOR, J. E. and SIEGMUND, D. (2004). Strong control, conservative point estimation and simultaneous conservative consistency of false discovery rates: A unified approach. *J. R. Stat. Soc. Ser. B Stat. Methodol.* **66** 187–205. MR2035766
[20] STUTE, W. (1982). The oscillation behavior of empirical processes. *Ann. Probab.* **10** 86–107. MR0637378
[21] TOBLER, W. R. (1979). Cellular geography. In *Philosophy in Geography* (S. Gale and G. Olsson, eds.) 379–386. Reidel Publishing Company, Dordrecht.
[22] WINKLER, G. (2003). *Image Analysis, Random Fields and Markov Chain Monte Carlo Methods*, 2nd ed. Springer, Berlin. MR1950762
[23] WU, W. B. (2005). Nonlinear system theory: Another look at dependence. *Proc. Natl. Acad. Sci. USA* **102** 14150–14154. MR2172215
[24] WU, W. B. (2006). Empirical processes of stationary sequences. *Statist. Sinica.* To appear.

DEPARTMENT OF STATISTICS
UNIVERSITY OF CHICAGO
5734 S. UNIVERSITY AVENUE
CHICAGO, ILLINOIS 60637
USA
E-MAIL: wbwu@galton.uchicago.edu